\documentclass[10pt]{amsart}

\usepackage{amsmath,amsthm,amsfonts,amscd,amssymb}
\usepackage{hyperref,wasysym}
\usepackage{color}

\newtheorem{thm}{Theorem}[section]
\newtheorem{cor}[thm]{Corollary}
\newtheorem{lem}[thm]{Lemma}
\newtheorem{prop}[thm]{Proposition}

\theoremstyle{definition}

\theoremstyle{remark}
\newtheorem{rem}{Remark}[section]
\newtheorem{ex}{Example}[section]

\begin{document}

\title{On sparse geometry of numbers}
\author{Lenny Fukshansky}\thanks{Fukshansky was partially supported by the Simons Foundation grant \#519058}
\author{Pavel Guerzhoy}
\author{Stefan K\"uhnlein}

\address{Department of Mathematics, 850 Columbia Avenue, Claremont McKenna College, Claremont, CA 91711}
\email{lenny@cmc.edu}
\address{Department of Mathematics, University of Hawaii, 2565 McCarthy Mall, Honolulu, HI, 96822-2273}
\email{pavel@math.hawaii.edu}
\address{Institut f\"ur Algebra und Geometrie, Fakult\"at f\"ur Mathematik, KIT,  FRG-76128 Karlsruhe}
\email{stefan.kuehnlein@kit.edu}

\subjclass[2010]{Primary: 11H06, 52C07, 11G05}
\keywords{lattices, sparse vectors, virtually rectangular lattices, Siegel's lemma, elliptic curve, $j$-invariant, isogeny, modular curve, geodesics}

\begin{abstract}
Let $L$ be a lattice of full rank in $n$-dimensional real space. A vector in $L$ is called $i$-sparse if it has no more than $i$ nonzero coordinates. We define the $i$-th successive sparsity level of $L$, $s_i(L)$, to be the minimal $s$ so that $L$ has $s$ linearly independent $i$-sparse vectors, then $s_i(L) \leq n$ for each $1 \leq i \leq n$. We investigate sufficient conditions for $s_i(L)$ to be smaller than $n$ and obtain explicit bounds on the sup-norms of the corresponding linearly independent sparse vectors in~$L$. These results can be viewed as a partial sparse analogue of Minkowski's successive minima theorem. We then use this result to study virtually rectangular lattices, establishing conditions for the lattice to be virtually rectangular and determining the index of a rectangular sublattice. We further investigate the $2$-dimensional situation, showing that virtually rectangular lattices in the plane correspond to elliptic curves isogenous to those with real $j$-invariant. We also identify planar virtually rectangular lattices in terms of a natural rationality condition of the geodesics on the modular curve carrying the corresponding points.
\end{abstract}

\maketitle

\def\A{{\mathcal A}}
\def\B{{\mathcal B}}
\def\C{{\mathcal C}}
\def\D{{\mathcal D}}
\def\F{{\mathcal F}}
\def\x{{\mathcal H}}
\def\I{{\mathcal I}}
\def\J{{\mathcal J}}
\def\K{{\mathcal K}}
\def\L{{\mathcal L}}
\def\M{{\mathcal M}}
\def\O{{\mathcal O}}
\def\R{{\mathcal R}}
\def\s{{\mathcal S}}
\def\V{{\mathcal V}}
\def\W{{\mathcal W}}
\def\X{{\mathcal X}}
\def\Y{{\mathcal Y}}
\def\H{{\mathcal H}}
\def\OO{{\mathcal O}}
\def\BB{{\mathbb B}}
\def\cee{{\mathbb C}}
\def\pee{{\mathbb P}}
\def\que{{\mathbb Q}}
\def\real{{\mathbb R}}
\def\zed{{\mathbb Z}}
\def\hyp{{\mathbb H}}
\def\aa{{\mathfrak a}}
\def\qbar{{\overline{\mathbb Q}}}
\def\eps{{\varepsilon}}
\def\ahat{{\hat \alpha}}
\def\bhat{{\hat \beta}}
\def\gt{{\tilde \gamma}}
\def\h{{\tfrac12}}
\def\be{{\boldsymbol e}}
\def\bei{{\boldsymbol e_i}}
\def\bff{{\boldsymbol f}}
\def\ba{{\boldsymbol a}}
\def\bb{{\boldsymbol b}}
\def\bc{{\boldsymbol c}}
\def\bm{{\boldsymbol m}}
\def\bk{{\boldsymbol k}}
\def\bi{{\boldsymbol i}}
\def\bl{{\boldsymbol l}}
\def\bq{{\boldsymbol q}}
\def\bu{{\boldsymbol u}}
\def\bt{{\boldsymbol t}}
\def\bs{{\boldsymbol s}}
\def\bv{{\boldsymbol v}}
\def\bw{{\boldsymbol w}}
\def\bx{{\boldsymbol x}}
\def\bX{{\boldsymbol X}}
\def\bz{{\boldsymbol z}}
\def\bwy{{\boldsymbol y}}
\def\bY{{\boldsymbol Y}}
\def\bL{{\boldsymbol L}}
\def\baa{{\boldsymbol\alpha}}
\def\bbb{{\boldsymbol\beta}}
\def\bet{{\boldsymbol\eta}}
\def\bxi{{\boldsymbol\xi}}
\def\bo{{\boldsymbol 0}}
\def\bol{{\boldkey 1}_L}
\def\ep{\varepsilon}
\def\p{\boldsymbol\varphi}
\def\q{\boldsymbol\psi}
\def\rank{\operatorname{rank}}
\def\aut{\operatorname{Aut}}
\def\lcm{\operatorname{lcm}}
\def\sgn{\operatorname{sgn}}
\def\spn{\operatorname{span}}
\def\md{\operatorname{mod}}
\def\Norm{\operatorname{Norm}}
\def\dim{\operatorname{dim}}
\def\det{\operatorname{det}}
\def\Vol{\operatorname{Vol}}
\def\rk{\operatorname{rk}}
\def\Gal{\operatorname{Gal}}
\def\WR{\operatorname{WR}}
\def\WO{\operatorname{WO}}
\def\GL{\operatorname{GL}}
\def\SL{\operatorname{SL}}
\def\PGL{\operatorname{PGL}}
\def\PSL{\operatorname{PSL}}
\def\adj{\operatorname{adj}}

\section{Introduction}
\label{intro}

Let $n \geq 2$ be an integer. For each $\bx \in \real^n$, we write
$$\|\bx\| = \left( \sum_{i=1}^n x_i^2 \right)^{1/2},\ |\bx| = \max_{1 \leq i \leq n} |x_i|$$
for the usual Euclidean norm and sup-norm on~$\real^n$, respectively. We also define the {\it $0$-norm} on~$\real^n$:
$$\|\bx\|_0 := \sum_{i=1}^n x_i^0,$$
where we use the convention that $0^0=0$. The $0$-norm counts the number of nonzero coordinates of a vector, which we refer to as the {\it sparsity level} of this vector; if sparsity level of some vector is no larger than $m$, we say that this vector is {\it $m$-sparse}. Sparsity has been actively investigated in the context of compressed sensing, which is a signal recovery paradigm based on the idea that most signals are sparse and can therefore be reconstructed from a small number of linear measurements~\cite{cs}. More recently, the sparsity phenomenon has also been studied in discrete mathematics and discrete geometry, in particular in the context of lattices~\cite{FNS}, \cite{aliev1}, \cite{aliev2}. In this paper, we want to take a first stab at a systematic approach to what we see as a ``sparse analogue" of the classical geometry of numbers.

Let $A = (a_{ij}) \in \GL_n(\real)$, and define
$$|A| := \max_{1 \leq i,j \leq n} |a_{ij}|.$$
Let $L = A\zed^n \subset \real^n$, then $L$ is a lattice of full rank with basis matrix~$A$. The minimal norm of $L$ is defined as
$$|L| := \min \left\{ \|\bx\| : \bx \in L \setminus \{\bo\} \right\}.$$
Previous research has focused on sparsity of integer representations of lattice vectors, i.e. on representing a vector $\bx \in L$ as $\bx = A\bwy$ with $\bwy \in \zed^n$ being as sparse as possible. In this paper, we will focus on the sparsity of the lattice vectors themselves. Specifically, we define the {\it successive sparsity levels} $s_1,\dots,s_n$ of the lattice $L$ to be
\begin{eqnarray*}
s_i(L) & := & \min \big\{ s : \exists\ i \text{ linearly independent vectors } \bx_1,\dots,\bx_i \in L \\
& & \text{ with } \|\bx_1\|_0, \dots, \|\bx_i\|_0 \leq s \big\}.
\end{eqnarray*}
Then $1 \leq s_1 \leq \dots \leq s_n \leq n$. Given a lattice $L$, what can be said about its successive sparsity levels? Further, assuming we know that some $s_{\ell} \leq k$, can we find $\ell$ such $k$-sparse vectors in~$L$? To answer these questions, we need some more notation. 

For every nonzero vector $\bx \in L$ define
$$d(\bx) := \dim_{\que} \spn_{\que} \left\{ x_{1},\dots,x_{n} \right\}$$
to be its {\it rational dimension}. If $A$ is an $n \times n$ real matrix with row vectors $\ba_i$ for $1 \leq i \leq n$, then for each subset $I \subseteq [n] := \{1,\dots,n\}$ we define $d_I(A) := \sum_{i \in I} d(\ba_i)$. We write $d(A)$ for $d_{[n]}(A)$, and define $d(L) = d(A)$, where $A$ is any basis matrix for $L$. Indeed, this definition does not depend on the choice of a basis matrix: if $A$ and $B$ are two basis matrices for $L$, then $B = AU$ for some $U \in \GL_n(\zed)$, and so each row vector $\bb_i$ of $B$ is of the form $\bb_i = \ba_i U$ for the corresponding row vector $\ba_i$ of $A$, which implies that $d(\bb_i) = d(\ba_i)$. More generally, $d(A) = d(AU)$ holds for any $U \in \GL_n(\que)$. Notice that $d(L) \geq n$. We refer to $d(L)$ as the rational dimension of $L$: the smaller $d(L)$ is the ``closer" $L$ is to being {\it rational}, meaning $L \subset \que^n$. Recall that $L$ is {\it integral} if $\|\bx\|^2 \in \zed$ for any $\bx \in L$, and $L$ is {\it arithmetic} if it is a scalar multiple of an integral lattice, so rational lattices are arithmetic. Certainly $d(L)=n$ for all rational lattices, but there also exist non-rational arithmetic lattices for which $d(L) = n$, for instance
$$L = \begin{pmatrix} \sqrt{2} & 2 \sqrt{2} \\ \sqrt{2} & 3 \sqrt{2} \end{pmatrix} \zed^2$$
is one such example. On the other hand, there exist arithmetic lattices with rational dimension $> n$, for instance
$$\begin{pmatrix} 1 & \sqrt{2} \\ \sqrt{2} & -1 \end{pmatrix} \zed^2,\ \begin{pmatrix} 1 & 1 & 1 \\ 1 & \sqrt{2} & 0 \\ \sqrt{2} & -1 & 0 \end{pmatrix} \zed^3,$$
and similar examples can be constructed in every dimension. Non-arithmetic lattices can also have rational dimension $> n$, for instance the planar lattice
\begin{equation}
\label{L12}
\Lambda_1 = \begin{pmatrix} 1 & \sqrt{3} \\ 0 & 1 \end{pmatrix} \zed^2
\end{equation}
with $d(\Lambda_1) = 3$, as well as $= n$, for instance the planar lattice
$$\Lambda_2 = \begin{pmatrix} \pi & 2 \pi \\ 2 & 1 \end{pmatrix} \zed^2$$
with $d(\Lambda_2) = 2$.
\smallskip

We will also define two ``measures of irrationality" of vectors in $L$ and of $L$ itself. First, if $\bx \in L$ has $d(\bx) = k$, then it can be written as
\begin{equation}
\label{alpha_rep}
\bx = \sum_{i=1}^k \alpha_i \bff_i,
\end{equation}
where $\bff_1,\dots,\bff_k$ are integer vectors with relatively prime coordinates. Notice that this decomposition is unique only if $k=1$, for example
\begin{eqnarray*}
\left( 1, \sqrt{2}, -2 \sqrt{3} \right) & = & 1 \cdot \left( 1, 0, 0 \right) + \sqrt{2} \cdot \left( 0,1,0 \right) - 2 \sqrt{3} \cdot \left( 0,0,1 \right) \\
& = & 1 \cdot \left( 1, 0, 0 \right) + \left( \frac{\sqrt{2}}{2} - \sqrt{3} \right) \cdot \left( 0,1,1 \right) + \left( \frac{\sqrt{2}}{2} + \sqrt{3} \right) \cdot \left( 0,1,-1 \right).
\end{eqnarray*}
Then define
$$\nu(\bx) := \left\{ \begin{array}{ll}
|\alpha_1| & \mbox{if $k=1$,} \\
0 & \mbox{if $k > 1$}.
\end{array}
\right.$$
For our basis matrix $A$, we define
$$\nu(A) := \prod_{i=1}^n \nu(\ba_i),$$
and for the lattice $L = A\zed^n$, we let $\nu(L) = \nu(A)$. This definition does not depend on the choice of the basis matrix~$A$. Clearly, there are many lattices for which $\nu(L) = 0$, for example $\nu(\Lambda_1) = 0$, where $\Lambda_1$ is as in~\eqref{L12}. In fact, it is not difficult to show that $\nu(L) > 0$ if and only if $d(L) = n$ (see Lemma~\ref{d_L} below). 

Second, let $\left< L \right>$ be the additive abelian group generated by the entries of vectors of $L$, and suppose that $\left< L \right>$ has rank $k \geq 1$. Fix a basis $\baa = (\alpha_1,\dots,\alpha_k)$ for $\left< L \right>$ over $\zed$, then every $\bx \in L$ has a representation of the form~\eqref{alpha_rep} with vectors $\bff_1,\dots,\bff_k \in \zed^n$. Define a map $\Phi_{\baa} : L \to \real^{nk}$ by
\begin{equation}
\label{phi_def}
\Phi_{\baa} \left( \sum_{i=1}^k \alpha_i \bff_i \right) = \begin{pmatrix} \bff_1 \\ \vdots \\ \bff_k \end{pmatrix}.
\end{equation}
The map $\Phi_{\baa}$ is additive, and hence extends to a map $\real \otimes L = \real^n \to \real^{nk}$. We can then pull back the sup-norm $| \cdot |$ on $\real^{nk}$ to $\real^n$ under $\Phi_{\baa}$ by defining $|\bx|_{\Phi_{\baa}} := |\Phi_{\baa}(\bx)|$, this way obtaining a norm $| \cdot |_{\Phi_{\baa}}$ on $\real^n$, which can then be compared to the sup-norm on~$\real^n$. Specifically, we can define
\begin{equation}
\label{mu_def}
\mu(\baa) := \sup \left\{ \frac{|\bx|_{\Phi_{\baa}}}{|\bx|} : \bx \in \real^n \setminus \{\bo\}  \right\}.
\end{equation}
We can now state our first result.

\begin{thm} \label{main1} Let $A \in \GL_n(\real)$ and let $L = A\zed^n$. Fix a basis $\baa$ for $\left< L \right>$ as above and let $\mu(\baa)$ be as given in~\eqref{mu_def}. Let $1 \leq k < n$ and suppose that there exists a subset $I \subset [n]$ of $n-k$ distinct indices such that $d_I(A) < n$. Let $\ell = n-d_I(A)$. Then $s_\ell(L) \leq k$, and there exist $\ell$ linearly independent vectors $\bx_1,\dots,\bx_{\ell} \in L$ with $\|\bx_i\|_0 \leq k$ and
\begin{equation}
\label{sm1}
\prod_{i=1}^{\ell} |\bx_i| \leq n^{n - d_I(A)/2} |A|^{n} \mu(\baa)^{d_I(A)}.
\end{equation}
\end{thm}

\noindent
This theorem can be viewed as a ``sparse" partial analogue of Minkowski's successive minima theorem. Indeed, if we know that $s_{\ell}(L) \leq k$, we can define the {\it $k$-sparse successive minima} $\lambda_1(L,k) \leq \dots \leq \lambda_{\ell}(L,k)$ with respect to sup-norm to be
$$\lambda_i(L,k) := \min \left\{ t \in \real_{>0} : \exists\ \text{ lin. ind. } \bx_1,\dots,\bx_i \in L \text{ with } \|\bx_j\|_0 \leq k, |\bx_j| \leq t \right\},$$
so the usual successive minima are $\lambda_i(L) := \lambda_i(L,n)$. Then~\eqref{sm1} is an upper bound on the product of these $k$-sparse successive minima. We prove Theorem~\ref{main1} in Section~\ref{sparsity}. Our main tool here is the celebrated Siegel's lemma, which is known to be sharp with respect to the exponent (we state it below as Theorem~\ref{siegel-lem}). We comment on the quality of the bound~\eqref{sm1} at the end of Section~\ref{sparsity}.  Unfortunately, the upper bound of~\eqref{sm1} depends on the choice of the basis for $L$ and for $\left< L \right>$. We can alleviate this dependence for lattices $L$ with rational dimension $d(L) = n$. We need some more notation.

Let us say that a lattice is {\it rectangular} if it has an orthogonal basis. Following~\cite{kuhnlein}, we will say that a lattice is {\it virtually rectangular} if it contains a rectangular sublattice of finite index. Two lattices $L, L' \subset \real^n$ are called {\it isometric} if there exists a real orthogonal matrix $U$ such that $L' = UL$; on the other hand, $L,L'$ are called {\it similar} if there exists a real orthogonal matrix $U$ and a positive real number $\beta$ such that $L' = \beta UL$. In other words, similarity as a linear map is a composition of an isometry and a dilation. It is easy to notice that the virtually rectangular property is preserved under isometry (and under similarity), however the sparsity levels, rational dimension and the irrationality measure $\nu$ are not necessarily preserved under isometry (they are preserved under dilation). In Section~\ref{rectangular} we give the following characterization of virtually rectangular lattices, using the invariants we have just introduced.

\begin{thm} \label{virt_rect} Let $L \subset \real^n$ be a lattice of full rank. The following three statements are equivalent:
\begin{enumerate}
\item $d(L) = n$,
\item $\nu(L) > 0$,
\item $s_1(L) = \dots = s_n(L) = 1$.
\end{enumerate}
Further, a full-rank lattice $L' \subset \real^n$ is virtually rectangular if and only if it is isometric to some lattice $L$ satisfying the three equivalent conditions above.
\end{thm}

In Section~\ref{rectangular} we also prove the following effective result and show it to be optimal (Example~\ref{ex_main2}).

\begin{thm} \label{main2} Let $A \in \GL_n(\real)$ be such that the lattice $L = A\zed^n$ satisfies the equivalent conditions of Theorem~\ref{virt_rect}. Then $L$ contains a rectangular sublattice $M$ with a basis of $1$-sparse vectors so that
\begin{equation}
\label{L_M_index}
[L : M] = \left( \frac{\det (L)}{\nu(L)} \right)^{n-1}.
\end{equation}
More generally, if $L' \subset \real^n$ be a virtually rectangular lattice, then there exists a rectangular sublattice $M'$ of $L'$ such that $[L':M'] = \left( \frac{\det (L)}{\nu(L)} \right)^{n-1}$, where $L$ is a lattice isometric to $L'$ which satisfies the equivalent conditions of Theorem~\ref{virt_rect}. 
\end{thm} 
\smallskip

In the 2-dimensional case our results imply a certain property of elliptic curves over~$\cee$. An elliptic curve $E$ can be realized as a complex torus $\cee/\Lambda$ for a planar lattice $\Lambda$, called the {\it period lattice} of this curve. Given two elliptic curves, $E$ and $E'$ a morphism $\phi : E \to E'$ between them such that $\phi(0) = 0$ and $\phi(E) \neq \{ 0 \}$ is called an {\it isogeny}. It is a remarkable fact that an isogeny is always surjective and has a finite kernel, the order of which is called the {\it degree} of the isogeny, denoted $\delta(E/E')$. If an isogeny $E \to E'$ exists, then there also exists the dual isogeny $E' \to E$ of the same degree such that their composition is simply the multiplication-by-$\delta(E/E')$ map, and hence the curves are called {\it isogenous}: this is an equivalence relation.  An injective isogeny is called an {\it isomorphism}, and the set of isomorphism classes of elliptic curves over~$\cee$ is parameterized by
\begin{equation}
\label{D_set}
\D := \{ \tau = a + bi \in \cee : -1/2 < a \leq 1/2, b \geq 0, |\tau| \geq 1 \} \setminus \{ e^{i\theta} : \pi/2 < \theta < 2\pi/3 \}
\end{equation}
in the following way. For each $\tau = a+bi \in \D$ we can define a lattice
\begin{equation}
\label{gtau}
\Gamma_\tau = \zed + \zed\tau,
\end{equation}
which can be thought of as $\begin{pmatrix} 1 & a \\ 0 & b \end{pmatrix} \zed^2$ in $\real^2$. Then every elliptic curve is isomorphic to an elliptic curve $E_{\tau}$ with period lattice $\Gamma_{\tau}$ for some $\tau \in \D$. In fact, it is more natural to identify the set of isomorphic classes of elliptic curves with the quotient space of the upper half-plane under the action of $\SL_2(\zed)$ by linear fractional transformations, where $\D$ is a fundamental domain for this action. There is a unique bijective holomorphic map $j : \D \to \cee$ taking $e^{2\pi i/3}$ to $0$ and $i$ to $1728$, called the {\it Klein $j$-function}, which is modular and gives the {\it $j$-invariant} $j(\tau)$ for each isomorphism class~$E_{\tau}$ of elliptic curves. In~\cite{pasha} the relevant properties of the $j$-invariant are outlined, and in particular it is noted that for $\tau \in \D$, $j(\tau) \in \real$ if and only if $\tau$ belongs to the set
\begin{equation}
\label{bnd}
\left\{ 1/2 + it : t \in \real, t \geq \sqrt{3}/2 \right\} \cup \left\{ e^{i\theta} : \theta \in [\pi/3,\pi/2] \right\} \cup \left\{ it : t \in \real, t \geq 1 \right\},
\end{equation}
and $j$ maps the first of these three subsets bijectively onto the interval $(-\infty,0]$, the second onto $[0,1]$, and the third onto $[1,\infty)$ (see also Proposition on p. 160 of~\cite{kk} for an earlier appearance of this observation). 
With this notation, we can state the following result which we prove in Section~\ref{ell_sect}.

\begin{thm} \label{elliptic} Let $\tau = a+bi \in \D$ and let $E_{\tau}$ be the corresponding elliptic curve with the period lattice $\Gamma_{\tau}$ as above. The following statements are equivalent:
\begin{enumerate}

\item Either $a \in \que$ or there exists some $t \in \real$ such that $a-bt, a+b/t \in \que$,

\item $\Gamma_{\tau}$ is virtually rectangular,

\item $E_{\tau}$ is isogenous to an elliptic curve $E'$ with real $j$-invariant $\geq 1$,

\item $E_{\tau}$ is isogenous to an elliptic curve $E'$ with real $j$-invariant in $[0,1]$.

\end{enumerate}
If these equivalent conditions hold with $a \in \que$, then there exists such an isogeny $E' \to E_{\tau}$ with $\delta(E'/E_{\tau}) = $ the denominator of $a$. If the conditions hold with $a \notin \que$ and $t$ is any real number satisfying (1), then there exists such an isogeny $E' \to E_{\tau}$ with
\begin{equation}
\label{delta}
\delta(E'/E_{\tau}) = \frac{|b| vw (t^2+1)}{|t|},
\end{equation}
where $v, w > 0$ are denominators of the rational numbers $a-bt$ and $a+b/t$, respectively.
\end{thm}

\noindent
Our proof of this theorem uses Theorem~\ref{virt_rect}. In particular, condition (1) of Theorem~\ref{elliptic} is equivalent to condition (1) of Theorem~\ref{virt_rect} in this 2-dimensional situation. Further, \eqref{delta} is just a reformulation of~\eqref{L_M_index} in this case, since $\delta(E'/E_{\tau})$ is precisely the index of the rectangular period lattice of $E'$ as a sublattice in the virtually rectangular period lattice~$\Gamma_{\tau}$ of~$E_{\tau}$. In the case when $a = p/q \in \que$, $d(\Gamma_{\tau}) = 2$.
\smallskip

We will refer to elliptic curves satisfying the conditions of Theorem~\ref{elliptic} as {\it virtually rectangular}. The class of their period lattices includes all $\Gamma_{\tau}$ so that $j(\tau) \in \real$ (see~\eqref{bnd}). Further, this class includes all arithmetic planar lattices (see Lemma~2.5 of~\cite{kuhnlein}), which are $\Gamma_{\tau}$ for $\tau \in \D$ being a quadratic irrationality (see~\cite{pasha}): these are $\tau = a + bi$ with $a, b^2 \in \que$, which correspond precisely to elliptic curves with complex multiplication (CM). We will discuss this situation in more details in Section~\ref{ell_sect}, in particular proving that CM elliptic curves are the only ones whose period lattice contains non-parallel rectangular sublattices (Proposition~\ref{t_unique} and Corollary~\ref{CM}): in the CM case, there are infinitely many $t$ satisfying condition (1) of Theorem~\ref{elliptic} (each corresponding to a different rectangular sublattice), whereas for all other virtually rectangular elliptic curves such $t$ is essentially unique. Finally, in Section~\ref{modular} we will show that virtually rectangular lattices in the plane have intrinsic geometric meaning in terms of the corresponding points on the modular curve: they correspond precisely to the points that lie on geodesics closed at infinity (Theorem~\ref{geodesic}).
\bigskip

\section{Successive sparsity levels}
\label{sparsity}

In this section we prove Theorem~\ref{main1}. We start with a lemma on successive sparsity levels.

\begin{lem} \label{irr} Let $A \in \GL_n(\real)$ and let $L = A\zed^n$. Let $1 \leq k < n$ and suppose that there exists a subset $I \subset \{1,\dots,n\}$ of $n-k$ distinct indices such that $d_I(A) < n$. Let $\ell = n-d_I(A)$. Then $s_\ell(L) \leq k$.
\end{lem}

\proof
Let $I = \{ i_1,\dots,i_{n-k} \}$ for some $1 \leq i_1 < i_2 < \dots < i_{n-k} \leq n$, and let us write $d_j := d(\ba_{i_j})$ for each $1 \leq j \leq n-k$. Let $A_I$ be the $(n-k) \times n$ submatrix of $A$ consisting of the rows indexed by~$I$. We want to show that there exists a nonzero vector $\bx \in \zed^n$ such that $A_I \bx = \bo$. For each $1 \leq j \leq n-k$, let
$$V_j = \left\{ \bx \in \que^n : \ba_{i_j} \cdot \bx = 0 \right\},$$
then $\dim_{\que} V_j = n - d_j$. Further, let us prove that
$$\dim_{\que} \left( \bigcap_{j=1}^{n-k} V_j \right) \geq \ell = n - d_I(A).$$
We argue by induction on $n-k \geq 1$. If $n-k = 1$, then $I = \{i_1\}$ and so $d_I(A) = d_1$, in which case 
$$\dim_{\que} V_1 = n-d_1 = \ell.$$
Now assume the result for all $1 \leq n-k < m \leq n-1$, and let us prove it for $n-k=m$. Let
$$I' = \{ i_1,\dots,i_{m-1} \}, \text{ so } I = I' \cup \{i_m\},$$
then $d_I(A) = d(I') + d_m$. Let $V' = \bigcap _{j=1}^{m-1} V_j$ and $V = V' \cap V_m$. By induction hypothesis,
$$\dim_{\que} V' \geq n - d(I').$$
Since $d(I') < d_I(A) < n$, this implies that $\dim_{\que} V' > 0$, and so $V' \neq \{\bo\}$. Now, by a well-known identity in linear algebra,
\begin{eqnarray*}
\dim_{\que} V & = & \dim_{\que} V' + \dim_{\que} V_m - \dim_{\que} \spn_{\que} \{ V', V_m \} \\
& \geq & (n - d(I')) + (n - d_m) - \dim_{\que} \spn_{\que} \{ V', V_m \} \\
& \geq & n - d_I(A) = \ell,
\end{eqnarray*}
since $\spn_{\que} \{ V', V_m \} \subseteq \que^n$, and so $\dim_{\que} \spn_{\que} \{ V', V_m \} \leq n$.

This implies that $\dim_{\que} \left( \bigcap_{j=1}^{n-k} V_j \right) = \ell  > 0$, and so there exist $\ell$ nonzero linearly independent vectors $\bwy_1,\dots,\bwy_\ell \in \bigcap_{j=1}^{n-k} V_j$. These vectors are in $\que^n$ and satisfy the equation $A_I \bwy_i = \bo$. Multiplying $\bwy_1,\dots,\bwy_\ell$ by the least common denominator of their coordinates, we obtain linearly independent vectors $\bx_1,\dots,\bx_\ell \in \zed^n$ such that $A_I \bx_i = \bo$. This means that the vectors $A\bx_1,\dots,A\bx_\ell \in L$ have at least $n-k$ coordinates equal to~$0$. Since $\bx_1,\dots,\bx_\ell$ are linearly independent and $A$ is a nonsingular matrix, we must have $A\bx_1,\dots,A\bx_\ell$ linearly independent, and so $s_\ell(L) \leq k$.
\endproof

\begin{rem} Notice that
$$d_I(A) \geq \dim_{\que} A_I := \dim_{\que} \{ a_{i_1 1},\dots,a_{i_k n} \}.$$
The converse of Lemma~\ref{irr} is not true: if $s_1(L) = k$, there may not exist any $I \subset \{1,\dots,n\}$ of cardinality $n-k$ so that $d_I(A) < n$. Indeed, consider the example
\begin{equation}
\label{A_matrix}
A = \begin{pmatrix} 1 & \sqrt{3} & 2 \sqrt{3} \\ \sqrt{5} & \sqrt{3} & 2 \sqrt{3} \\ \sqrt{2} & \sqrt{3} & \sqrt{5} \end{pmatrix}
\end{equation}
and let $L = A \zed^3$. Then $s_1(L) = 1$, since
$$A \begin{pmatrix} 0 \\ 2 \\ -1 \end{pmatrix} = \begin{pmatrix} 0 \\ 0 \\ 2 \sqrt{3} - \sqrt{5} \end{pmatrix} \in L.$$
On the other hand, $d(\ba_1) = d(\ba_2) = 2$ and $d(\ba_3) = 3$. Thus for any $I$ of cardinality $3-1 = 2$, $d_I(A) \geq 4 > n = 3$. Further, even $\dim_{\que} A_I$ here is at least $3$. On the other hand, notice for comparison purposes that if a lattice $L = A\zed^n$ is virtually rectangular, then $\dim_{\que} (A^{\top} A) \leq n$ (see~\cite{kuhnlein}). 
\end{rem}

For each row vector $\ba_i$ of $A$ let $d_i = d(\ba_i)$. Then there exist $\que$-linearly independent real numbers $\alpha_{i1},\dots,\alpha_{id_i}$ such that
\begin{equation}
\label{faa}
\ba_i = \sum_{j=1}^{d_i} \alpha_{ij} \bff_{ij},
\end{equation}
where $\bff_{ij}$ are integer vectors with relatively prime coefficients for all $1 \leq i \leq n$, $1 \leq j \leq d_i$. Let $d = \sum_{i=1}^n d_i$ and let $F(A)$ be the $d \times n$ matrix with rows $\bff_{ij}$. 

\begin{ex}
For instance, in case of the matrix $A$ in~\eqref{A_matrix}, we have $d_1 = 2$, $d_2 = 2$, $d_3 = 3$, and
$$\ba_1 = \alpha_{11} \bff_{11} + \alpha_{12} \bff_{12},\ \ba_2 = \alpha_{21} \bff_{21} + \alpha_{22} \bff_{22},\ \ba_3 =  \alpha_{31} \bff_{31} +  \alpha_{32} \bff_{32} + \alpha_{33} \bff_{33},$$
where
$$\alpha_{11} = 1,\ \alpha_{12} = \alpha_{22} = \alpha_{32} = \sqrt{3},\ \alpha_{21} = \alpha_{33} = \sqrt{5},\ \alpha_{31} = \sqrt{2},$$
and
$$\bff_{11} = \bff_{21} = \bff_{31} = (1,0,0),\ \bff_{12} = \bff_{21} = (0,1,2),\ \bff_{32} = (0,1,0),\ \bff_{33} = (0,0,1).$$
Therefore $d = 2+2+3 = 7$ in this example, and the $7 \times 3$ matrix $F(A)$ is
$$F(A) = \begin{pmatrix} 1 & 0 & 0 \\ 0 & 1 & 2 \\ 1 & 0 & 0 \\ 0 & 1 & 2 \\ 1 & 0 & 0 \\ 0 & 1 & 0 \\ 0 & 0 & 1 \end{pmatrix}.$$
\end{ex}

Define
$$|F(A)| := \max_{1 \leq i \leq n} \max_{1 \leq j \leq d_i} |\bff_{ij}|,$$
where $|\bff_{ij}|$ is the sup-norm of the vector $\bff_{ij}$.

\begin{lem} \label{f_coeff} For a fixed choice of $\baa = \left( \alpha_{ij} : 1 \leq i \leq n, 1 \leq j \leq d_i \right)$ as above, we have
$$|F(A)| \leq \mu(\baa) |A|,$$
where $\mu(\baa)$ is as in~\eqref{mu_def}.
\end{lem} 

\proof
Notice that the rank of the additive group $\left< L \right>$ is equal to $d$, and $\baa$ is a basis for it. Using the notation of Section~\ref{intro}, specifically~\eqref{phi_def}, notice that for each $1 \leq i \leq n$,
$$\frac{|\ba_i|_{\Phi_{\baa}}}{|\ba_i|} \leq \mu(\baa).$$
On the other hand, $|\ba_i|_{\Phi_{\baa}} = |\Phi_{\baa}(\ba_i)| = \max_{1 \leq j \leq d_i} |\bff_{ij}|$, and hence
$$|F(A)| = \max_{1 \leq i \leq n} \max_{1 \leq j \leq d_i} |\bff_{ij}| \leq \mu(\baa) \max_{1 \leq i \leq n} |\ba_i| = \mu(\baa) |A|.$$
\endproof

Our next result relies heavily on the use of Siegel's lemma (Theorem~2 of~\cite{siegel}) and its adaptation to sup-norm using Fisher's inequality (equation (1.8) of~\cite{masser}). We state it here for the reader's convenience.

\begin{thm} [Siegel's lemma with matrix sup-norm]  \label{siegel-lem} Let $B$ be an $m \times n$ integer matrix of rank $m < n$. Then there exist $n-m$ linearly independent vectors $\bz_1,\dots,\bz_{n-m} \in \zed^n$ such that $B\bz_i = \bo$ for every $1 \leq i \leq n-m$ and
$$\prod_{i=1}^{n-m} |\bz_i| \leq \left( \sqrt{n} |B| \right)^m.$$
The exponent $m$ in this upper bound cannot in general be improved.
\end{thm}

\begin{lem} \label{succ_min} Let $A \in \GL_n(\real)$, $L = A\zed^n$ and let $\baa$ be a fixed basis for $\left< L \right>$. Let $1 \leq k < n$ and suppose that there exists a subset $I \subset \{1,\dots,n\}$ of $n-k$ distinct indices 
$$1 \leq i_1 < i_2 < \dots < i_{n-k}$$
such that $d_I(A) := \sum_{j=1}^{n-k} d_{i_j} < n$. Let $\ell = n-d_I(A)$. Then there exist $\ell$ linearly independent vectors $\bx_1,\dots,\bx_{\ell} \in L$ with $\|\bx_i\|_0 \leq k$ and
$$\prod_{i=1}^{\ell} |\bx_i| \leq n^{n - d_I(A)/2} |A|^{n} \mu(\baa)^{d_I(A)}.$$
\end{lem} 

\proof
Let $A_I$ be the $(n-k) \times n$ submatrix of $A$ consisting of the rows indexed by~$I$. Let $F(A)_I$ be the $d_I(A) \times n$ submatrix of $F(A)$ consisting of rows with $\bff_{i_l j}$ for $i_l \in I$ and $1 \leq j \leq d_{i_l}$. Notice that $A_I \bwy = \bo$ for some $\bwy \in \zed^n$ if and only if $F(A)_I \bwy = 0$. Using the notation in the proof of Lemma~\ref{irr}, we have
$$V = \left( \bigcap_{j=1}^{n-k} V_j \right) = \left\{ \bwy \in \que^n : A_I \bwy = \bo \right\} = \left\{ \bwy \in \que^n : F(A)_I \bwy = \bo \right\},$$
which is an $\ell$-dimensional subspace of $\que^n$. By Theorem~\ref{siegel-lem}, there exist $\ell$ linearly independent vectors $\bwy_1,\dots,\bwy_{\ell} \in V \cap \zed^n$ such that
\begin{equation}
\label{siegel_lem}
\prod_{j=1}^{\ell} |\bwy_j| \leq \left( \sqrt{n} |F(A)| \right)^{d_I(A)}.
\end{equation}
For each $1 \leq j \leq \ell$, define $\bx_j = A\bwy_j$. Since $A$ is a nonsingular matrix, $\bx_1,\dots,\bx_{\ell}$ are nonzero linearly independent vectors in $L$ which are at least $k$-sparse. Now notice that
\begin{equation}
\label{x-to-y}
|\bx_j| = |A\bwy_j| \leq n |A| |\bwy_j|,
\end{equation}
and so
$$\prod_{i=1}^{\ell} |\bx_i| \leq n^{\ell} |A|^{\ell} \prod_{j=1}^{\ell} |\bwy_j|.$$
Combining this observation with~\eqref{siegel_lem} and Lemma~\ref{f_coeff} completes the proof.
\endproof

Theorem~\ref{main1} now follows by combining Lemmas~\ref{irr} and~\ref{succ_min}. Notice that the exponent $d_I(A)$ in the upper bound of~\eqref{siegel_lem} is sharp due to the optimality of the bound in Theorem~\ref{siegel-lem}.  On the other hand, dependence of the bound of~\eqref{x-to-y} on $|A|$ also cannot be improved in general. This suggests that the dependence of the bound of Theorem~\ref{main1} on $|A|$ and $\mu(\baa)$ has correct order of magnitude. 
\bigskip

\section{Virtually rectangular lattices}
\label{rectangular}

In this section we focus on virtually rectangular lattices. We start by presenting the proof of Theorem~\ref{virt_rect}, split into two parts.

\begin{lem} \label{d_L} Let $L \subset \real^n$ be a lattice of full rank. The following three statements are equivalent:
\begin{enumerate}

\item $d(L) = n$,

\item $\nu(L) > 0$,

\item $s_1(L) = \dots = s_n(L) = 1$.

\end{enumerate}

\end{lem}

\proof
Let $L = A\zed^n$, where $A$ is a basis matrix with rows $\ba_1,\dots,\ba_n$. We will prove that (1) is equivalent to (2) and that (1) is equivalent to (3). First assume that $d(L) = n$, then $d(\ba_i) = 1$ for each row $\ba_i$ of the basis matrix~$A$. This means that each $\ba_i = \alpha_i \bz_i$, where $\alpha_i \in \real \setminus \{0\}$ and $\bz_i \in \zed^n$ is a vector with relatively prime coordinates, and so $\nu(\ba_i) = |\alpha_i|$. Then
$$\nu(L) = \nu(A) = \prod_{i=1}^n |\alpha_i| > 0,$$
and so (1) implies (2). Further, this means that for any subset $I \subset [n]$ of cardinality $n-1$, the linear system $A_I \bwy = \bo$ has a nontrivial integer solution. Since such sets $I$ are of the form $I = [n] \setminus \{ i\}$, $1 \leq i \leq n$, we can index corresponding integer solutions by $\bwy_i$. Then each vector $A\bwy_i$ has only the $i$-th coordinate nonzero, and hence all of such vectors are linearly independent (they are multiples of the standard basis vectors). Therefore $s_1(L) = \dots = s_n(L) = 1$, and so (1) implies (3).

Next assume $\nu(L) > 0$, and let $A$ be a basis matrix for $L$. Then
$$\nu(L) = \nu(A) = \prod_{i=1}^n \nu(\ba_i) > 0,$$
which means that each $\ba_i$ is of the form $\ba_i = \alpha_i \bz_i$ for some $\alpha_i \in \real \setminus \{0\}$ and $\bz_i \in \zed^n$ a primitive vector. Hence $d(L) = n$, and so (2) implies (1).

Finally, suppose $s_1(L) = \dots = s_n(L) = 1$. Then there exist linearly independent vectors $a_1 \be_1,\dots,a_n \be_n \in L$, where $\be_1,\dots,\be_n$ are the standard basis vectors. Hence there exists $U \in \GL_n(\que)$ such that $AU$ is a nonsingular diagonal matrix, which implies $d(A) = d(AU) = n$. Thus (3) implies (1).
\endproof

\begin{lem} \label{rect} A lattice $L$ is virtually rectangular if and only if it is isometric to some lattice $L'$ with $s_1(L') = \dots = s_n(L') = 1$.
\end{lem}

\proof
Suppose that $L$ contains a rectangular sublattice $M$, and let $B$ be an orthogonal basis matrix for $M$. Then there exists a real orthogonal matrix $U$ such that $UB$ is a diagonal matrix. Let $M' = UM = UB\zed^n$ be a sublattice of the lattice $L' = UL$. Since $UB$ is diagonal, $M'$ has a basis consisting of scalar multiples of the standard basis vectors, and thus all successive sparsity levels of $L'$ are equal to~$1$.

Conversely, assume $L$ is isometric to some $L'$ with successive sparsity levels equal to~$1$, say $L = UL'$ for some orthogonal matrix $U$. Then $L'$ contains $n$ linearly independent vectors $\bx_1,\dots,\bx_n$ with $\|\bx_i\|_0 = 1$. These vectors must therefore be constant multiples of standard basis vectors. Let $M' = \spn_{\zed} \{ \bx_1,\dots,\bx_n \}$, then $M = UM'$ is a rectangular sublattice of~$L$.
\endproof

\noindent
Then Theorem~\ref{virt_rect} follows by combining Lemmas~\ref{d_L} and~\ref{rect}. Next we prove Theorem~\ref{main2}.

\proof[Proof of Theorem~\ref{main2}]
Since $d(L) = n$, we must have $d(\ba_i) = 1$ for each row vector $\ba_i$ of $A$. Then there must exist nonzero real numbers $\alpha_1,\dots,\alpha_n$ and primitive integer row vectors $\bff_1,\dots,\bff_n$ so that $\ba_i = \alpha_i \bff_i$ for each $1 \leq i \leq n$.  Hence the matrix $F(A)$ with row vectors $\bff_i$ is $n \times n$ and $A = \A F(A)$, where $\A$ is the diagonal matrix with diagonal entries $\alpha_1,\dots,\alpha_n$. Let $\adj(F(A))$ be the adjugate of $F(A)$, then $\adj(F(A))$ is an integer matrix in $\GL_n(\que)$, $\det(\adj(F(A))) = \det(F(A))^{n-1}$ and
$$F(A) \adj(F(A)) = \det (F(A)) I_n,$$
where $I_n$ is the $n \times n$ identity matrix. Then
$$B := A \adj(F(A)) = \A F(A) \adj(F(A)) = \det (F(A)) \A,$$
which is a diagonal matrix with diagonal entries $\alpha_1 \det(F(A)),\dots,\alpha_n \det(F(A))$. This implies that
$$\det(B) = \det(F(A))^n \det(\A) = \det(F(A))^n \prod_{i=1}^n \alpha_i.$$
On the other hand, since $\adj(F(A))$ is an integer matrix in $\GL_n(\que)$, the lattice $M := B \zed^n$ is a full-rank sublattice of $L$, which is rectangular with a basis of $1$-sparse vectors. Further,
\begin{eqnarray*}
[L : M] & = & \frac{\det(M)}{\det(L)} = \frac{|\det(B)|}{|\det(A)|} = | \det(\adj(F(A)) | \\
& = & | \det(F(A)) |^{n-1} = \left| \frac{\det(A)}{\det(\A)} \right|^{n-1} = \left( \frac{\det(L)}{\nu(L)} \right)^{n-1}.
\end{eqnarray*}
Now, if $L' \subset \real^n$ is any virtually rectangular lattice, then it is isometric to some lattice $L$ satisfying the equivalent conditions of Theorem~\ref{virt_rect}. This $L$ has a rectangular sublattice $M$ as we just constructed satisfying~\eqref{L_M_index}. Then $L'$ contains a rectangular sublattice $M'$ isometric to $M$ with $[L':M'] = [L:M]$. 
\endproof

\begin{ex} \label{ex_main2}
Theorem~\ref{main2} is optimal, i.e. the lattice $L = A\zed^n$ as in the statement of the theorem may not contain a rectangular sublattice $M$ with a basis of $1$-sparse vectors and smaller index than given in~\eqref{L_M_index}. Indeed, this is easily seen to be the case when
$$A =  \begin{pmatrix} 1 & 0 & \dots & 0 & 0 \\ 
0 & 1 & \dots & 0 & 0 \\
\vdots & \vdots & \ddots & \vdots & \vdots \\
0 & 0 & \dots & 1 & 0 \\
1 & 1 & \dots & 1 & d \end{pmatrix}$$
for some integer $d > 1$. Here $\nu(L) = 1$, $\det(L) = d$, and the smallest-index rectangular sublattice with a basis of $1$-sparse vectors is
$$M = \begin{pmatrix} d & 0 & \dots & 0 & 0 \\ 
0 & d & \dots & 0 & 0 \\
\vdots & \vdots & \ddots & \vdots & \vdots \\
0 & 0 & \dots & d & 0 \\
0 & 0 & \dots & 0 & d \end{pmatrix} \zed^n,$$
which has index $d^{n-1}$ in $L$.
\end{ex}

\begin{rem} \label{2-dim-pf} It may be instructive to separately consider the $2$-dimensional case of Theorem~\ref{main2}, where the computation becomes completely elementary. Let $L' \subset \real^2$ be a virtually rectangular lattice and let $L$ be a lattice isometric to $L'$ which satisfies the equivalent conditions of Theorem~\ref{virt_rect}. Then $L = A\zed^2$, where
$$A = \begin{pmatrix} \alpha_1 u_1 & \alpha_1 v_1 \\ \alpha_2 u_2 & \alpha_2 v_2 \end{pmatrix}$$
with $u_1,u_2,v_1,v_2$ relatively prime integers and $\alpha_1,\alpha_2$ real numbers. Then
$$\det(L) = \alpha_1 \alpha_2 |u_1v_2 - u_2v_1|,\ \nu(L) = | \alpha_1 \alpha_2 |.$$
Further, it is easy to see that the orthogonal vectors
$$\bz_1 = \begin{pmatrix} 0 \\ \alpha_2(u_1v_2 - u_2v_1) \end{pmatrix},\ \bz_2 = \begin{pmatrix} \alpha_1 (u_2v_1 - u_1v_2) \\ 0 \end{pmatrix}$$
are in $L$, and so $M = \spn_{\zed} \{ \bz_1,\bz_2 \}$ is a rectangular sublattice of $L$. Then
$$\det(M) = \left| \alpha_1 \alpha_2 (u_1v_2 - u_2v_1) (u_2v_1 - u_1v_2) \right| =  \frac{\det(L)^2}{| \alpha_1 \alpha_2 |}.$$
Let $M'$ be a sublattice of $L'$ isometric to $M$ in $L$, then
$$[L':M'] = [L:M] = \frac{\det(M)}{\det(L)} = \frac{\det(L)}{\nu(L)}.$$
\end{rem}
\bigskip

\section{Isogenies of elliptic curves}
\label{ell_sect}

In this section we prove Theorem~\ref{elliptic} and discuss some of its consequences. To start with, we state a technical lemma that will be of use to us: it is a combination of Lemma~5.3 and Proposition~5.4 of~\cite{pasha} (see also Proposition on p. 160 of~\cite{kk}, as mentioned in Section~\ref{intro}).

\begin{lem} \label{53} Let $\D$ be as in~\eqref{D_set}. For $\tau \in \D$, the value $j(\tau)$ is real if and only if $\tau$ belongs to the set described in~\eqref{bnd}. Further, $\Gamma_{\tau}$ is WR if and only if $j(\tau)$ is real and belongs to the interval $[0,1]$.
\end{lem}

\proof[Proof of Theorem~\ref{elliptic}]
First suppose that $a = \frac{p}{q} \in \que$, then $\Gamma_{\tau}$ contains orthogonal vectors
$$\begin{pmatrix} 1 \\ 0 \end{pmatrix},\ \begin{pmatrix} 0 \\ qb \end{pmatrix} = q \begin{pmatrix} a \\ b \end{pmatrix} - p \begin{pmatrix} 1 \\ 0 \end{pmatrix}.$$
These two vectors span a rectangular sublattice of $\Gamma_{\tau}$ of determinant $qb$, i.e. of index $q$, which in particular implies that $\Gamma_{\tau}$ is virtually rectangular. If $a \notin \que$, assume that there exists some $t \in \real$ such that $a-bt, a+b/t \in \que$. Define the lattice
$$L_t := \frac{1}{\sqrt{1+t^2}} \begin{pmatrix} 1 & a-bt \\ t & at+b \end{pmatrix} \zed^2,$$
then it is easy to see that $d(L_t) = 2$, and so it is virtually rectangular by Theorem~\ref{virt_rect}. Let $\theta = \arctan t$, then $\cos \theta = \frac{1}{\sqrt{1+t^2}}$ and $\sin \theta = \frac{t}{\sqrt{1+t^2}}$, meaning that
$$U_t = \frac{1}{\sqrt{1+t^2}} \begin{pmatrix} 1 & -t \\ t & 1 \end{pmatrix}$$
is an orthogonal matrix. Notice that $U_t \Gamma_{\tau} = L_t$, meaning that $\Gamma_{\tau}$ is isometric to $L_t$, hence it is also virtually rectangular. This shows that condition~(1) implies~(2).

Suppose now $\Gamma_{\tau}$ is virtually rectangular, then it contains a rectangular sublattice~$\Gamma'$. Let $E'$ be the elliptic curve (up to isomorphism) with period lattice $\Gamma'$ (up to similarity). We can then assume that 
\begin{equation}
\label{gprime}
\Gamma' = \begin{pmatrix} 1 & 0 \\ 0 & q \end{pmatrix} \zed^2,
\end{equation}
which means that $E' = E_{\tau'}$ for $\tau' = iq$. Hence $\tau'$ is in the third component set of~\eqref{bnd}, and so $j(\tau') \geq 1$ (see Lemma~\ref{53} above). Now, since the period lattice of $E'$ is a sublattice of the period lattice of $E_{\tau}$, there must exist an isogeny $E' \to E_{\tau}$ induced by the projection $\cee/\Gamma' \to \cee/\Gamma_{\tau}$. This shows that condition~(2) implies~(3).

An isogeny $E' \to E$ exists if and only if the period lattice for  $E'$ is (up to similarity) a sublattice of the period lattice for $E$. A planar lattice is called {\it well-rounded (WR)} if it has two linearly independent shortest vectors with respect to Euclidean norm. Now, the period lattice is rectangular if and only if the corresponding $j$-invariant is real and $\geq 1$  while the period lattice is WR if and only if the corresponding $j$-invariant is real and in the interval $[0,1]$ (Lemma~\ref{53}). Hence to prove that conditions (3) and (4) are equivalent it is sufficient to show that $\Gamma_{\tau}$ contains a rectangular sublattice if and only if it contains a WR sublattice. This is guaranteed by Lemma~2.1 of~\cite{kuhnlein}.

Next assume that $E_{\tau}$ is isogenous to some elliptic curve $E' = E_{\tau'}$ with real nonnegative $j$-invariant. Let $\Gamma' = \Gamma_{\tau'}$ be the period lattice of $E'$, so $j(\tau') \in \real_{\geq 0}$ and $\Gamma'$ is (up to similarity) a sublattice of $\Gamma$. If $j(\tau') \geq 1$, then Lemma~\ref{53} guarantees that $\tau' = iq$ for some $q \geq 1$, and so $\Gamma'$ is of the form~\eqref{gprime}, which is rectangular. If, on the other hand, $0 \leq j(\tau') < 1$, then Lemma~\ref{53} implies that the lattice $\Gamma'$ is WR. Now, Lemma~2.1 of~\cite{kuhnlein} asserts that a lattice has WR sublattices if and only if it is virtually rectangular. Hence we conclude in any case that $\Gamma_{\tau}$ is virtually rectangular. Therefore $\Gamma_{\tau}$ is isometric to some lattice $L$ with $d(L) = 2$, by Theorem~\ref{virt_rect}. Let
$$U(\theta) = \begin{pmatrix} \cos \theta & -\sin \theta \\ \sin \theta & \cos \theta \end{pmatrix}$$
for some angle $\theta$ be the corresponding isometry matrix, and let $\tau = a+bi$ so that $\Gamma_{\tau}$ is of the form~\eqref{gtau}. Then
$$L = U(\theta) \Gamma_{\tau} = \begin{pmatrix} \cos \theta & a \cos \theta - b \sin \theta \\ \sin \theta & a \sin \theta + b \cos \theta \end{pmatrix} \zed^2 = \frac{1}{\sqrt{1+t^2}} \begin{pmatrix} 1 & a-bt \\ t & at+b \end{pmatrix} \zed^2,$$
where $t = \tan \theta$. Since $d(L) = 2$, we must have $a - bt \in \que$ and $\frac{at+b}{t} = a + b/t \in \que$. This shows that condition~(3) implies~(1).
\medskip

Finally, assume that the equivalent conditions of Theorem~\ref{elliptic} hold. If $a = \frac{p}{q} \in \que$, then $\Gamma_{\tau}$ contains a rectangular sublattice $\Gamma'$ of index~$q$. Let $E'$ be the elliptic curve (up to isomorphism) corresponding to $\Gamma'$, then the degree of the isogeny $E' \to E_{\tau}$ is precisely this index $q$. If $a \notin \que$, then equivalent conditions of Theorem~\ref{elliptic} hold with some $t \in \real$. Then the period lattice $\Gamma_{\tau}$ of the curve $E_{\tau}$ is virtually rectangular and isometric to the lattice $L_t = A_t \zed^2$ with
$$A_t = \frac{1}{\sqrt{1+t^2}} \begin{pmatrix} 1 & a-bt \\ t & at+b \end{pmatrix} =  \begin{pmatrix} \frac{1}{\sqrt{1+t^2}}  & \frac{1}{\sqrt{1+t^2}} (a-bt) \\ \frac{t}{\sqrt{1+t^2}}  & \frac{t}{\sqrt{1+t^2}} \left( a+\frac{b}{t} \right) \end{pmatrix},$$
and $d(L'_t) = 2$. Since $a-bt, a+b/t \in \que$, we can write
$$a-bt =\frac{u}{v},\ a+\frac{b}{t} = \frac{q}{w}$$
with $u,v,q,w \in \zed$ and $v,w > 0$. Then, repeating the argument of Remark~\ref{2-dim-pf} for this specific situation,
$$u - (a-bt)v = 0,\ q - \left( a + \frac{b}{t} \right) w = 0,$$
and so the vectors
$$\frac{t}{\sqrt{1+t^2}} \begin{pmatrix} 0 \\ u - \left( a + \frac{b}{t} \right) v \end{pmatrix},\ \frac{1}{\sqrt{1+t^2}} \begin{pmatrix} q - (a-bt) w \\ 0 \end{pmatrix}$$
are in $L_t$. These two vectors span a rectangular sublattice $R_t$ of $L_t$, whose determinant is
$$\det R_t = \left| \frac{t}{1+t^2} \left( u - \left( a + \frac{b}{t} \right)  v \right) (q - (a-bt) w) \right| = \frac{b^2 vw (t^2+1)}{|t|}.$$
Let $\Gamma'$ be a rectangular sublattice of $\Gamma_t$ isometric to $R_t$, then
$$[\Gamma_t : \Gamma'] = [L_t : R_t] = \frac{\det R_t}{\det L_t} = \frac{|b| vw (t^2+1)}{|t|}.$$
Now, let $E'$ be the elliptic curve (up to isomorphism) corresponding to $\Gamma'$ (up to similarity), then the degree of the isogeny $E' \to E_{\tau}$ is
$$\delta(E'/E_{\tau}) = [\Gamma_{\tau} : \Gamma'] = \frac{|b| vw (t^2+1)}{|t|}.$$
This completes the proof.
\endproof
\smallskip

Notice that if a lattice $\Gamma_{\tau}$ for $\tau = a+bi \in \D$ contains a rectangular sublattice, then it contains infinitely many non-similar rectangular sublattices: these can be obtained for instance by multiplying the original rectangular sublattice by matrices of the form $\begin{pmatrix} l & 0 \\ 0 & m \end{pmatrix}$ for relatively prime integers $l,m$. However all of these sublattices are parallel to each other, meaning that they are spanned by parallel pairs of orthogonal basis vectors. Can $\Gamma_{\tau}$ have non-parallel rectangular sublattices? This condition is equivalent to saying that there are multiple ways to rotate $\Gamma_{\tau}$ so that some sublattice will have an orthogonal basis along the coordinate axes. Since the parameter $t$ of Theorem~\ref{elliptic} is the tangent of the angle of rotation, this can be possible if and only if there exist distinct $t_1,t_2 \in \real$ satisfying condition (1) of Theorem~\ref{elliptic} so that $t_1 \neq -1/t_2$ ($t$ and $-1/t$ correspond to rotations resulting in the same lattice). This turns out to be possible if and only if $\tau$ is a quadratic irrationality, in which case the corresponding elliptic curve~$E_{\tau}$ is said to be a curve with {\it complex multiplication (CM)}: this is precisely the situation when the endomorphism ring of $E_{\tau}$ is larger than~$\zed$ (specifically, an order in the imaginary quadratic field~$\que(\tau)$; see Corollary III.9.4 of~\cite{silverman}). We now prove that for such $\tau$ there are infinitely many different real numbers $t$ satisfying condition (1) of Theorem~\ref{elliptic}.

\begin{prop} \label{t_unique} With notation as in Theorem~\ref{elliptic}, suppose that there exist $t_1,t_2 \in \real$ satisfying condition (1). Define
\begin{align}
\alpha_1 & := a-bt_1 \in \que \label{eq:1}\\
\beta_1 & := a+b/t_1 \in \que \label{eq:2}\\
\alpha_2 & := a-bt _2  \in \que \label{eq:3}\\
\beta_2 & := a+b/t_2 \in \que  \label{eq:4}
\end{align}
Assume also that $t_1 \neq t_2, -1/t_2$. Then $t_1^2,t_2^2 \in \que$ and $b^2 \in \que$, meaning that $\tau = a+bi$ is a quadratic irrationality, and hence $E_\tau$ is a CM elliptic curve. 
\end{prop}

\proof
Subtract (\ref{eq:1}) from  (\ref{eq:3}) and (\ref{eq:2}) from  (\ref{eq:4}) to obtain
\begin{align}
b(t_1-t_2) & = \alpha_1-\alpha_2  \label{eq:2.1}\\
b\frac{t_1-t_2}{t_1t_2} &= \beta_1-\beta_2. \label{eq:2.2}
\end{align}
Divide (\ref{eq:2.1}) by (\ref{eq:2.2}) (note that $t_1 \neq t_2$ implies that $\beta_1 \neq \beta _2$) to conclude that
\[
t_1t_2 = \frac{\alpha_1-\alpha_2}{\beta_1-\beta_2} \in \que.
\]
Multiply    (\ref{eq:4}) by $t_1t_2$ and take into the account $bt_1 = a - \alpha_1$ (from (\ref{eq:1})) to obtain
\[
a(t_1t_2 + 1) = \alpha_1 +\beta_2 t_1t_2
\]
and conclude that $a \in \que$ since $t_1t_2 \neq -1$ by assumption, and the quantities $\alpha_1,\beta_2 \in \que$, and $t_1t_2\in \que$ are already known to be rational. Since $a \in \que$, we conclude from (\ref{eq:1}) and (\ref{eq:3}) that $bt_1, bt_2 \in \que$, and therefore their ratio (note that $t_1t_2 \neq 0$ because $\tau \notin \real$) is rational: $t_1/t_2 \in \que$. Multiplication (and division) by the rational quantity $t_1t_2$ now allows us to conclude that $t_1^2, t_2^2 \in \que$. Finally, since $bt_1 = a - \alpha_1 \in \que$, we square it to conclude that $b^2t_1^2 \in \que$, and divide by $t_1^2$ to obtain that $b^2 \in \que$ as claimed. 
\endproof

\noindent
Proposition~\ref{t_unique} asserts essential uniqueness of the real number $t$ satisfying condition (1) of Theorem~\ref{elliptic} in the generic (non-CM) situation. In the case when CM occurs, there is an infinite family of such~$t$.

\begin{cor} \label{CM} With notation as in Proposition~\ref{t_unique}, assume $\tau = a + bi$ with $a,b^2 \in \que$. Then every $t$ satisfying condition (1) of Theorem~\ref{elliptic} is of the form 
\begin{equation}
\label{q_cond}
qb \text{ or } q/b \text{ for some } q \in \que.
\end{equation}
Conversely, for every rational $q$, $t = qb$ satisfies condition (1) of Theorem~\ref{elliptic}.
\end{cor}

\proof
Assume $t_1, t_2$ are two different values of $t$ satisfying condition (1) of Theorem~\ref{elliptic}. From the proof of Proposition~\ref{t_unique}, we know that $t_1/t_2 \in \que$, so $t_1 = qt_2$ for some $q \in \que$. Then $t_1$ satisfies~\eqref{q_cond} if and only if $t_2$ does, so it is enough to show that there exists a $t$ of the form~\eqref{q_cond} satisfying condition (1) of Theorem~\ref{elliptic}. Indeed, for any $q \in \que$, take $t = qb$, then
$$a - bt = a - b^2q \in \que,\ a + b/t = a + 1/q \in \que.$$
On the other hand, take~$t = q/b$, then
$$a - bt = a - q \in \que,\ a + b/t = a + b^2/q \in \que.$$
This completes the proof.
\endproof

Let us now provide an interpretation of this result for elliptic curves. For every $N > 1$ there is a symmetric polynomial $F_N(X,Y)=F_N(Y,X)$ in two variables with integer coefficients, which has the following property: two elliptic curves $E_1$ and $E_2$ with corresponding $j$-invariants $j_1$ and $j_2$ are isogenous with an isogeny of degree $N$ if and only if $F_N(j_1,j_2)=0$. The polynomial $F_N(X,Y)$ is commonly referred to as {\it $N$-th modular polynomial}.
Our Theorem~\ref{elliptic} now implies that an elliptic curve $E$ with a non-real $j$-invariant $j(E)$ is virtually rectangular if and only if for some $N$ the polynomial $F_N(X,j(E))$, now monic with complex coefficients, has a real root. Another observation that is not difficult to prove is that the curve $E_{\tau}$ is virtually rectangular if and only if $E_{-\bar{\tau}}$ is virtually rectangular, and if this is the case then the curves $E_{\tau}$ and $E_{-\bar{\tau}}$ are isogenous (i.e., any one of the corresponding lattices is similar to a sublattice of the other one).

Further, consider an elliptic curve $E$ over $\cee$ with a non-real $j$-invariant $j(E)$ which is not isomorphic to an elliptic curve over~$\real$. Assume that $E$ is nevertheless isogenous to an elliptic curve $E'$ over $\real$ with $j$-invariant $j(E')$. This implies that $j(E)$ and $j(E')$ are algebraically dependent over~$\que$. Then the upper bound on the inequality~\eqref{delta} on Theorem~\ref{elliptic} gives a bound on the degree of the field extension $[\que(j(E),j(E')): \que(j(E))]$.

\bigskip

\section{Planar virtually rectangular lattices on the modular curve}
\label{modular}

Our goal here is to present a geometric interpretation which justifies the consideration of virtually rectangular lattices as very natural objects. Specifically, we look at how the points corresponding to the virtually rectangular lattices are positioned in the moduli space of all lattices. 

WR lattices can be clearly seen in the fundamental domain (see Lemma~\ref{53} above), as can be the rectangular ones: they correspond to the set  $\{ it : t \in \real, t \geq 1 \}$. Virtually WR lattices (those that contain a finite-index WR sublattice) are the same as virtually rectangular in~$\real^2$ (see Lemma~2.1 of~\cite{kuhnlein}), however they are not as easily identified in the fundamental domain. Indeed, for a rational $a \in \que$ with $-1/2 < a \leq 1/2$, the lattice corresponding to every point $\tau=a+ib \in \D$ in the fundamental domain  is virtually rectangular. While this set is already too large to be contained on any continuous path in the upper half-plane, there are many more points in the fundamental domain corresponding to virtually rectangular lattices. Surprisingly, the picture becomes much clearer if we look at these points as points on the modular curve instead. 

Let 
$$\mathfrak H = \{\tau =x+iy \in \cee \ \vert \ y=\Im(\tau)>0\}$$
be the complex upper half-plane. It comes with the Poincar\'e metric 
$$ds^2=y^{-2}(dx^2+dy^2),$$
which is invariant under the action of $\PGL_2(\real)^+$ (and is uniquely defined by that property up to a constant multiplication). Let $Y = \PSL_2(\zed) \backslash {\mathfrak H}$. The points of $Y$ classify elliptic curves up to isomorphism over $\cee$ and correspond to (orientation preserving) similarity classes of lattices $\Gamma_\tau = \langle 1,\tau \rangle_{\zed}$ in the plane.

Since $\PSL_2(\zed) \subset \PGL_2(\real)$, the space $Y$ inherits the metric from $\mathfrak H$, in particular, geodesics on $Y$ are precisely the images of the geodesics on ${\mathfrak H}$ under the natural projection $\pi: {\mathfrak H} \rightarrow Y$ (recall that geodesics for a given metric are the paths of shortest length). The modular curve is the compact Riemann surface $X = \PSL_2(\zed) \backslash \overline{\mathfrak H}$, where $ \overline{\mathfrak H} = {\mathfrak H} \cup \que \cup\{\infty\}$.
We then have $X=Y\cup \infty$, and the point at infinity $\infty$ does not correspond to any lattice. 

Geodesics on $\mathfrak H$ (for the metric $ds^2$) are the vertical lines together with the semicircles orthogonal to the real axis (see e.g. \cite[Proposition 4.5.5]{Cohen-Stromberg}, \cite[Lemma 1.4.1]{Miyake}).  While all geodesics on $\mathfrak H$ are of infinite length, under the map $\pi$, some geodesics become closed, while others are still of infinite length. We say that a geodesic on $\mathfrak H$ {\it passes through $\infty$} if it is either a vertical line or a semicircle which meets the real line at a rational point, and we say that a geodesic is {\it closed at $\infty$} if it is either vertical with a rational $x$-coordinate or both ends of the semicircle meet the real line at rational points. We apply the same terminology to the images of these geodesics under the projection $\pi$. 
This terminology is not standard. For example, while for any two points in $Y$ there exists exactly one geodesic which passes through these two points, there are infinitely many geodesics which pass through $\infty$ and any given point on $Y$; furthermore, the geodesics which are closed at $\infty$ are never closed in $Y$. However, one may possibly argue that, for example, a semicircle which meets the real line at two rational points is really closed at $\infty$ because both of its ends map to $\infty$ under $\pi$.   

\begin{thm} \label{geodesic} A point on $Y$ corresponds to a virtually rectangular lattice if and only if this point belongs to a geodesic that is closed at $\infty$.
\end{thm}

\proof
Let $p \in Y$, and assume that the corresponding lattice is virtually rectangular. 
Thus $p=\pi(\tau)$ with $\tau=a+bi \in  \mathfrak H$, and the  lattice $\Gamma_\tau$ is virtually rectangular. Then the lattice $\Gamma_\tau$ has two orthogonal vectors, say, $ A\tau + B$ and $C\tau +D $ where $A,B,C,D$ are integers and $AD-BC \neq 0$. The orthogonality condition can be written as
\begin{equation} 
\label{eq:quadratic}
(a^2+b^2)AC+a(AD+BC) +BD =0.
\end{equation}
If $A=0$, then $BC\neq 0$, and $\tau$ belongs to the vertical line $x=-D/C$. If $C=0$, then $AD \neq 0$, and  $\tau$ belongs to the vertical line $x=-B/A$. Finally, if $AC \neq 0$ then (\ref{eq:quadratic}) becomes an equation of the semicircle
$$\left(a+\frac{AD+BC}{2AC} \right)^2 + b^2 = \frac{(AD-BC)^2}{4A^2C^2}$$
satisfied by $(a,b)$ with a rational radius of $|(AD-BC)/2AC|$ and the center on real line at $x=-(AD+BC)/2AC \in \que$. Thus in either case $\pi(\tau)$ belongs to a geodesic that is closed at $\infty$.

Conversely, assume that $\pi(\tau)$ belongs to a geodesic which is closed at $\infty$. If this geodesic is the image under $\pi$ of a vertical line in $\mathfrak H$ with a rational $x$-coordinate, then $\tau=a + bi$ with $a \in \que$, and the lattice $\Gamma_\tau$ is virtually rectangular. Otherwise, $\tau$ belongs to a semicircle which meets the real line at rational points, say $\alpha$ and $\beta$. Pick $\sigma \in \SL_2(\zed)$ such that $\sigma(\alpha)=\infty$. Then $\sigma(\beta) \in \que$, and $\sigma$ takes the semicircle to the vertical line with a rational $x$-coordinate of $\sigma(\beta)$ (M\"obius transformations preserve the Poincar\'e metric, therefore take geodesics to geodesics). Thus the lattice $\langle 1, \sigma(\tau) \rangle_{\zed}$ is virtually rectangular, and this lattice is similar to $\Gamma_\tau$.
We thus conclude that  $\Gamma_\tau$ is virtually rectangular.
\endproof
\bigskip

\noindent
{\bf Acknowledgement:} We wish to thank the anonymous referees for their helpful suggestions, which improved the quality of presentation.
\bigskip

\bibliographystyle{plain}  

\end{document}